\documentclass[12pt]{article}

\usepackage[latin1]{inputenc}
\usepackage{vmargin}
\usepackage[francais]{babel}
\usepackage[T1]{fontenc}
\usepackage[tbtags]{amsmath}
\usepackage{amsthm,amssymb}
\usepackage{pifont}
\usepackage[all]{xy}

\setmarginsrb{2cm}{1.5cm}{2cm}{1.5cm}{0.5cm}{0.5cm}{0.5cm}{1.2cm}

\DeclareFontFamily{T1}{rsfs}{}
\DeclareFontShape{T1}{rsfs}{m}{n}{ <-> rsfs10 }{}

\newtheorem{theo}{Théorème}[subsection]
\newtheorem{lemme}[theo]{Lemme}
\newtheorem{prop}[theo]{Proposition}

\newenvironment{preuve}{
\noindent \textbf{Démonstration.}}{\hfill \(\square\) }

\def\geq{\geqslant}

\def\pa#1{\left(#1\right)}

\def\acco#1{\left\{ #1 \right\}}

\def\D{\mathbb{D}}

\def\Z{\mathbb{Z}}

\def\J{\mathcal{J}}
\def\O{\mathcal{O}}
\def\S{\mathcal{S}}

\def\hom{{\text{Hom}}}
\def\ext{{\text{Ext}}}
\def\hhom{\text{\fontfamily{rsfs}\selectfont H}om}
\def\eext{\text{\fontfamily{rsfs}\selectfont E}xt}

\def\Fil{{\text{Fil}}}
\def\Frac{\text{Frac}\,}

\def\calA{\mathcal{A}}
\def\calM{\mathcal{M}}
\def\calN{\mathcal{N}}
\def\calX{\mathcal{X}}
\def\calF{\mathcal{F}}
\def\G{\mathcal{G}}
\def\H{\mathcal{H}}

\def\Gr{\text{Gr}}
\def\Mod{\text{Mod}}
\def\Mods{(\Mod / S)}
\def\pModS{'(\Mod / S)}
\def\ModS1{(\Mod / S_1)}
\def\Mr{\protect\underline{\mathcal{M}}^r}

\def\spec{\text{Spec\,}}
\def\DP{\text{DP}}
\def\pr{\text{pr}}
\def\id{{\text{id}}}
\def\im{{\text{im}}}

\def\cris{\text{cris}}
\def\CRIS{\text{CRIS}}
\def\syncris{\text{syn-cris}}
\def\SYNCRIS{\text{SYN-CRIS}}
\def\syn{\text{syn}}
\def\SYN{\text{SYN}}

\def\Ga{\mathbb{G}_{\text a}}
\def\Gm{\mathbb{G}_{\text m}}

\def\Phi{\phi}

\def\toinj{\hookrightarrow}

\makeatletter
\renewcommand{\everyentry@}{\vphantom{A_{[]}^{[]}}}
\makeatother

\title{Dualité de Cartier et modules de Breuil}
\author{Xavier Caruso}
\date{Septembre 2005}

\begin{document}

\maketitle

\renewcommand{\abstractname}{Abstract}
\begin{abstract}
Let $\O_K$ be a complete discrete valuation ring. Denote by $K$ its 
fractions field et by $k$ its residue field. Assume that $k$ is of 
characteristic $p > 0$ and perfect. In \cite{breuil}, Breuil gives
an anti-equivalence between the category of finite flat $\O_K$-group 
schemes killed by a power of $p$ and a category of linear algebra 
objects which is called $\Mods$. The aim of this article is to make
explicit the Cartier duality on the category $\Mods$.
\end{abstract}

\medskip

\renewcommand{\abstractname}{Résumé}
\begin{abstract}
Soit $K$ un corps complet pour une valuation discrète, de 
caractéristique nulle, et dont le corps résiduel est supposé parfait de 
caractéristique $p > 0$. On appelle $\O_K$ son anneau des entiers. Dans 
\cite{breuil}, Breuil exhibe une anti-équivalence de catégories entre
la catégorie des $\O_K$-schémas en groupes commutatifs finis, plats et 
tués par une puissance de $p$ et une certaine catégorie d'objets 
d'algèbre linéaire qu'il note $\Mods$. Le but de cet article est 
d'expliciter la dualité de Cartier sur la catégorie $\Mods$ \emph{via} 
l'équivalence précédente.
\end{abstract}

\tableofcontents

\vspace{1cm}

Tout au long de cet article, on considère $k$ un corps parfait 
de caractéristique $p > 0$. On note $W$ l'anneau des vecteurs de Witt à 
coefficients dans $k$, $\sigma : W \to W$ l'opérateur de Frobenius sur 
$W$ et $K_0 = \Frac W$. Soient $K$ une extension finie totalement 
ramifiée de $K_0$ et $\O_K$ son anneau des entiers. On fixe $\pi$ une 
uniformisante de $K$ et on note $E(u)$ le polynôme minimal de $\pi$ sur 
$K_0$.

\bigskip

Dans tout ce qui suit, on appelle \emph{$\O_K$-groupe} un 
$\O_K$-schéma en groupes commutatifs fini, plat et annulé par une 
puissance de $p$. Dans \cite{breuil}, Breuil construit une 
anti-équivalence de catégories entre la catégorie des $\O_K$-groupes 
(resp. la catégorie des groupes $p$-divisibles sur $\O_K$) et une 
catégorie d'objets d'algèbre linéaire qu'il note $\Mods$ (resp. une
catégorie d'objets appelés \emph{modules fortement divisibles}) dont la 
définition est rappelée en \ref{sec:defcat}.

Le but de cet article est d'expliciter purement en terme d'algèbre 
linéaire le foncteur déduit de la dualité de Cartier sur la catégorie 
$\Mods$ et sur la catégorie des modules fortement divisibles.

\bigskip

La première section est destinée au rappel des constructions et des 
résultats de 
\cite{breuil}. En particulier, on définit les catégories sus-mentionnées 
ainsi que le foncteur qui réalise l'anti-équivalence précédente. Dans la 
deuxième section, nous construisons une dualité sur les catégories 
d'objets d'algèbre linéaire et finalement nous prouvons dans la 
troisième section que cette dualité correspond \emph{via} le foncteur de 
Breuil à la dualité de Cartier sur les schémas en groupes.

\bigskip

Cet article est une version courte du chapitre V de la thèse de l'auteur 
(\cite{caruso-these}). Dans \emph{loc. cit.}, la construction de la 
dualité est étendue aux catégories $\Mr$ plus générales introduites par 
Breuil dans \cite{breuil-invent}. Ces dernières sont munies de foncteurs 
vers la catégorie des $\Z_p$-représentations galoisiennes, et on prouve 
(encore dans \cite{caruso-these}) certaines compatibilités entre la 
dualité et ces foncteurs.

Nous invitons le lecteur désireux d'avoir de nombreux compléments
à se reporter à cette référence.

\section{Rappel sur la classification de Breuil}

Dans cette section, on se contente de rappeler les résultats 
principaux de \cite{breuil}.

\subsection{Les catégories d'objets d'algèbre linéaire}
\label{sec:defcat}
Définissons l'anneau $S$ comme le complété $p$-adique de l'enveloppe à 
puissances divisées de $W[u]$ par 
rapport à l'idéal principal engendré par $E(u)$. On note $\Fil^1 S 
\subset S$ le complété $p$-adique de l'idéal engendré par les 
$\frac{E(u)^i}{i!}$ pour $i \geq 1$. On munit en outre $S$ d'un 
opérateur de Frobenius $\phi$ défini comme l'unique application $\phi : 
S \to S$ continue, $\sigma$-semi-linéaire et vérifiant $\phi(u) = u^p$.
On vérifie que $\phi(\Fil^1 S) \subset pS$, ce qui permet de définir 
$\phi_1 = \frac{\phi} p _{| \Fil^1 S}$. Finalement, on pose $c = \phi_1 
(E(u))$, c'est une unité de $S$.

\medskip

On définit à présent la catégorie $\pModS$. Ses objets sont la donnée 
d'un $S$-module $\calM$, d'un sous-module $\Fil^1 \calM \subset \calM$ 
tel que $\Fil^1 S \, \calM \subset \Fil^1 \calM$ et d'une application 
$\phi$-linéaire $\phi_1 : \Fil^1 \calM \to \calM$ vérifiant la condition :
$$\phi_1 (sx) = \frac 1 c \phi_1(s) \phi_1 (E(u) x)$$
pour tout $s \in \Fil^1 S$ et tout $x \in \calM$. Les morphismes de
$\pModS$ sont les applications $S$-linéaires respectant toutes les 
structures additionnelles. On a une notion de suite exacte dans cette 
catégorie : une suite est dite exacte si elle l'est en tant que suite 
de $S$-modules et si, en outre, la suite déduite sur les $\Fil^1$ est 
aussi exacte (comme suite de $S$-modules).

La catégorie $\ModS1$ est la sous-catégorie pleine 
de $\pModS$ formée des objets annulés par $p$ et tels que $\phi_1 
(\Fil^1 \calM)$ engendre $\calM$ en tant que $S$-module. Finalement,
la catégorie $\Mods$ est la plus petite sous-catégorie (pleine) de 
$\pModS$ contenant les objets de $\ModS1$ et stable par extension.

\bigskip

Un \emph{module fortement divisible} est par définition un objet de 
$\pModS$ pour lequel les trois conditions suivantes sont 
satisfaites :
\begin{dinglist}{43}
\item $\calM$ est libre de rang fini sur $S$,
\item $\calM / \Fil^1 \calM$ est sans $p$-torsion,
\item $\phi_1(\Fil^1 \calM)$ engendre $\calM$ en tant que $S$-module.
\end{dinglist}

\medskip

Un résultat important (et que l'on aura à manipuler par la suite)
concernant les modules fortement divisibles est donné par le lemme qui 
suit :

\begin{lemme}
\label{lem:adapt}
Soit $\calM$ un module fortement divisible. Il existe une base $(e_1, 
\ldots, e_d)$ de $\calM$ et des entiers $n_1, \ldots, n_d$ égaux à $0$ 
ou $1$ tels que :
$$\Fil^1 \calM = \bigoplus_{i=1}^d \Fil^{n_i} S \, e_i.$$
\end{lemme}

\begin{preuve}
C'est exactement le lemme 2.1.1.9 de \cite{breuil}.
\end{preuve}

\bigskip

\noindent
Une base vérifiant la condition de lemme précédent est appelée 
\emph{base adaptée (à la filtration)} de $\calM$. Ainsi le lemme affirme 
que tout module fortement divisible admet une base adaptée.

\subsection{L'anti-équivalence de catégories}
Dans cette partie, on donne la construction du foncteur $\Mod$ qui
s'avère être celui qui réalise l'anti-équivalence de catégories entre la
catégorie $\Mods$ et la catégorie des $\O_K$-groupes. On commence par
plusieurs rappels sur les topologies critalline et syntomique et sur les
principaux faisceaux pour ces topologies.

\subsubsection{Sites cristallin et syntomique}

On rappelle qu'un morphisme de schéma $X \to Y$ est dit 
\emph{syntomique} s'il est plat, localement de présentation finie
et s'il se factorise localement en une immersion fermée régulière
dans un $X$-schéma lisse. Les morphismes syntomiques sont stables par 
composition et changement de base.
Si $X$ est un schéma, on définit le gros (resp. le petit) site 
syntomique $X_\SYN$ (resp. $X_\syn$) comme la catégorie des $X$-schémas
(resp. des $X$-schémas syntomiques) munie de la topologie syntomique : 
une famille de morphismes $f_i : U_i \to U$ est un recouvrement si 
chacun des $f_i$ est syntomique et si topologiquement $\dot U = \bigcup 
f_i(\dot U_i)$.

\medskip

Soit $\Upsilon$ un schéma muni de puissances divisées 
et sur lequel $p$ est localement nilpotent. Si $X \to \Upsilon$ est 
tel que les puissances divisées sur $\Upsilon$ s'étendent à $X$, on
définit le site syntomique-cristallin (ou simplement cristallin) associé 
au morphisme $X \to \Upsilon$ de la façon suivante. La catégorie 
sous-jacente au site est l'ensemble des quadruplets $(U, T, i, \delta)$ 
tels que :
\begin{dinglist}{43}
\item $U$ est un schéma défini sur $X$,
\item $T$ est un schéma défini sur $\Upsilon$ sur lequel $p$ est 
localement nilpotent,
\item $i : U \toinj T$ est une immersion fermée définie sur $\Upsilon$,
\item $\delta$ est une structure d'idéal à puissances divisées sur 
l'idéal de $\O_T$ définissant l'immersion $i$, compatible aux 
puissances divisées sur $\Upsilon$.
\end{dinglist}
Par la suite, on écrira abusivement $(U,T)$ à la place de $(U, T, i, 
\delta)$. Un morphisme entre $(U,T)$ et $(U',T')$ est la donnée de deux
applications $U \to U'$ définie sur $X$ et $T \to T'$ définie sur 
$\Upsilon$ qui commutent aux $i$. Une famille $(U_i, T_i) \to (U,T)$ est 
un recouvrement si chacun des morphismes $U_i \to U$ est syntomique, si 
tous les diagrammes :
$$\xymatrix @C=50pt {
U_i \ar@{^(->}[r] \ar[d] & T_i \ar[d] \\
U \ar@{^(->}[r] & T }$$
sont cartésiens et finalement si topologiquement $\dot T = \bigcup_i 
f_i(\dot T_i)$ (où $f_i$ désigne le morphisme induit $T_i \to T$).

\medskip

Si on note $\widetilde{X_\SYN}$ et $(\widetilde {X/\Upsilon})_\SYNCRIS$ 
les catégories de faisceaux abéliens sur les deux sites précédents, on a 
un morphisme de topoï :
$$w : (\widetilde {X/\Upsilon})_\SYNCRIS \to \widetilde{X_\SYN}$$
donné par le couple de foncteurs adjoints $(w^\star, w_\star)$ définis 
par les formules suivantes :
\begin{eqnarray*}
w^\star \calF (U,T) & = & \calF(U) \\
w_\star \calF (U) & = & H^0 ((U/\Upsilon)_\SYNCRIS, 
\calF_{|(U/\Upsilon)_\SYNCRIS} ).
\end{eqnarray*}
On vérifie directement que $w_\star \circ w^\star = \id$, d'où on déduit 
que $w^\star$ est pleinement fidèle.

\bigskip

On dispose en outre de faisceaux importants sur les sites précédents. 
Sur le site syntomique, on montre que le préfaisceau $U \mapsto 
\Gamma(U, \O_U)$ est un faisceau que l'on appelle le \emph{faisceau 
structural}. 
Sur le site cristallin, on définit les faisceaux $\O_{X / 
\Upsilon}$ et $\Ga$ par les formules :
$$\O_{X / \Upsilon} (U,T) = \Gamma(T, \O_T) \qquad \text{et} \qquad
\Ga(U,T) = \Gamma(U, \O_U).$$
On a un morphisme naturel et surjectif $\O_{X / \Upsilon} \to \Ga$. On 
note $\J_{X / \Upsilon}$ son noyau, ce qui donne naissance 
tautologiquement à une suite exacte dans $(\widetilde 
{X/\Upsilon})_\SYNCRIS$ :
\begin{equation}
\label{5:eq:exactcris}
\xymatrix {
0 \ar[r] & \J_{X / \Upsilon} \ar[r] & \O_{X / \Upsilon} \ar[r] & \Ga 
\ar[r] & 0 }.
\end{equation}

\subsubsection{Les faisceaux $\O_n^\cris$ et $\J_n^\cris$}
\label{5:subsec:synlocal}

À partir de maintenant, on fixe un entier $n$, on note $S_n$ la 
réduction modulo $p^n$ de l'anneau $S$ introduit en \ref{sec:defcat} et 
$E_n$ le schéma $\spec(S_n)$. On pose également $T_n = 
\spec(\O_K/p^n)$. Les deux schémas précédents sont munis de puissances 
divisées (définies respectivement sur les idéaux $(E(u))$ et $(p)$) et 
on dispose d'un épaississement $T_n \toinj E_n$. On s'intéresse 
désormais plus particulièrement au cas $X = T_n$ et $\Upsilon = E_n$.

\bigskip

On note $\O_n$ le 
faisceau structural sur $(T_n/E_n)_\SYNCRIS$ défini par $\O_n = 
w_\star \Ga$. On définit pareillement $\O_n^\cris = w_\star 
\O_{T_n/E_n}$ et $\J_n^\cris = w_\star \J_{T_n/E_n}$. Il s'agit de 
faisceaux sur le gros site syntomique que l'on sait décrire localement 
sur la restriction au petit site syntomique.

\medskip

Soit $U$ un schéma syntomique sur $T_n$. Étale-localement, c'est le 
morphisme de schémas associé au morphisme d'anneaux $\O_K/p^n \to A$ 
avec :
$$A = \frac{\O_K/p^n[X_1, \ldots, X_s]}{(f_1, \ldots, f_t)}$$
où $X_1, \ldots, X_s$ sont des indéterminées et $(f_1, \ldots, f_t)$ une 
suite transversalement régulière relativement à $\O_K/p^n$. Posons pour 
tout entier $i$ :
$$A^i = \frac{\O_K/p^n[X_0^{1/p^i}, X_1^{1/p^i}, \ldots, 
X_s^{1/p^i}]}{(X_0-\pi, f_1, \ldots, f_t)}$$
et $A^\infty = \varinjlim A_i$ (pour les morphismes de transition 
évidents). Notons $W_n$ l'anneau des vecteurs de Witt de longueur $n$ à 
coefficients dans $k$ et $\phi$ le Frobenius sur cet anneau. Posons :
$$W_n^\cris (A^\infty) = W_n(A^\infty/p A^\infty) \otimes_{W_n,(\phi^n)} 
W_n[u].$$
On dispose d'une surjection $s : W_n^\cris (A^\infty) \to A^\infty$ qui 
envoie $u$ sur $X_0$ et $(a_0, \ldots, a_{n-1}) \in W_n(A^\infty/p 
A^\infty)$ sur $\hat a_0^{p^n} + p \hat a_1^{p^{n-1}} + \cdots + p^{n-1} 
\hat a_{n-1}^p$ où $\hat a_i \in A^\infty$ désigne un relevé de $a_i$.
On note $W_n^{\cris, \DP} (A^\infty)$ l'enveloppe à puissances divisées
de $W_n^\cris (A^\infty)$ par rapport au noyau de $s$ (et compatibles
avec les puissances divisées sur l'idéal $(p)$). La surjection $s$ se 
prolonge en une application $W_n^{\cris, \DP} (A^\infty) \to A^\infty$ 
que l'on note encore $s$.

\begin{lemme}
Avec les notations précédentes, il existe un isomorphisme canonique :
$$\varinjlim_i \O_n^\cris (A^i) \to W_n^{\cris, \DP} (A^\infty)$$
faisant commuter le diagramme suivant :
$$\xymatrix @C=50pt {
\varinjlim_i \O_n^\cris (A^i) \ar[d]^-{\sim} \ar[r] & \varinjlim_i \O_n 
(A^i) \ar@{=}[d] \\
W_n^{\cris, \DP} (A^\infty) \ar[r]^-{s} & A^\infty }$$
où la flèche du haut est obtenue en appliquant $w_\star$ au morphisme 
de faisceaux $\O_{T_n/E_n} \to \Ga$.
\end{lemme}

\begin{preuve}
Voir preuve du lemme 2.3.2 de \cite{breuil}.
\end{preuve}

\bigskip

On déduit directement du lemme précédent l'exactitude de la suite :
\begin{equation}
\label{5:eq:exactsyn}
\xymatrix {
0 \ar[r] & \J_n^\cris \ar[r] & \O_n^\cris \ar[r] & \O_n \ar[r] & 0 }
\end{equation}
obtenue en appliquant le foncteur $w_\star$ à la suite exacte 
(\ref{5:eq:exactcris}). Ce même lemme assure également que le faisceau 
$\O_n^\cris$ est plat sur $S_n$ tandis que $\J_n^\cris$ l'est sur $W_n$
(voir\footnote{Dans ce lemme, seule la platitude sur $W_n$ est annoncée
mais la platitude du $S_n$ est également vraie en reprenant les
arguments de la proposition 2.1.2.1 de \cite{breuil-duke}.} à nouveau le 
lemme 2.3.2 de \cite{breuil}). Ceci permet de vérifier que le 
Frobenius induit un morphisme $\phi : \J_{n+1}^\cris \to p 
\O_{n+1}^\cris$ qui s'annule sur $p^n \J_{n+1}^\cris$ et donc conduit à 
une flèche $\J_{n+1}^\cris/p^n \to p \O_{n+1}^\cris$. D'autre part, 
si $i : T_n \toinj T_{n+1}$ désigne l'épaississement évident, la même 
platitude fournit l'identification $\J_{n+1}^\cris/p^n \simeq 
i_\star  \J_n^\cris$ et montre que la multiplication par $p$ induit un
isomorphisme $i_\star \O_n^\cris \simeq p \O_{n+1}^\cris$. Il existe par
suite un unique morphisme $\phi_1$ qui fait commuter le diagramme 
suivant :
$$\xymatrix @C=50pt {
\J_{n+1}^\cris \ar[d]_-{\sim} \ar[r]^-{\phi} & p \O_{n+1}^\cris \\
i_\star \J_n^\cris \ar[r]^-{\phi_1} & i_\star \O_n^\cris 
\ar[u]^-{\sim}_-{p} }$$
Intuitivement, il faut penser à $\phi_1$ comme au quotient $\frac \phi 
p$, ou encore comme à un inverse du Verschiebung.

\subsubsection{Le foncteur $\Mod$}
\label{subsec:defmod}

Soit $\G$ un $\O_K$-groupe.
Pour $m = n$ et $m = n+1$, on note $\G_m = \G \times_{\spec(\O_K)} T_m$. 
Ces schémas définissent des faisceaux sur les sites $(T_n)_\syn$ (resp. 
$(T_n/E_n)_\SYNCRIS$) et $(T_{n+1})_\syn$ que l'on note encore $\G_n$ et 
$\G_{n+1}$. On appelle encore $i$ l'épaississement $T_n \toinj T_{n+1}$.
L'objet $\Mod(\G)$ est défini par (voir paragraphe 4.2.1 de 
\cite{breuil}) :
\begin{eqnarray*}
\Mod(\G) & = & \hom (\G_n, \O_n^\cris) = \hom (\G_{n+1}, i_\star 
\O_n^\cris) \\
\Fil^1 \Mod(\G) & = & \hom (\G_n, \J_n^\cris) = \hom (\G_{n+1}, i_\star 
\J_n^\cris) 
\end{eqnarray*}
et $\phi_1 : \Fil^1 \Mod(\G) \to \Mod(\G)$ est la flèche induite par 
$\phi_1 : i_\star \J_n^\cris \to i_\star \O_n^\cris$. Notons que les
$\hom$ précédents sont tous calculés dans la catégorie des 
\emph{faisceaux abéliens}. Dans la suite, il en sera toujours ainsi.

\bigskip

\noindent
{\it Remarques.} L'égalité $\hom (\G_n, \O_n^\cris) = \hom (\G_{n+1}, 
i_\star \O_n^\cris)$ résulte du fait que ces deux termes s'identifient
au même sous-ensemble de $\O_n^\cris(\G_n)$ puisque $\G_n$ (resp. 
$\G_{n+1})$ est syntomique sur $T_n$ (resp. $T_{n+1}$) (voir proposition 
2.2.2 de \cite{breuil}).

La définition de \cite{breuil} n'est pas exactement la même que celle 
que l'on vient de donner. En effet, dans \emph{loc. cit.}, il est 
question des faisceaux $\O_\infty^\cris$ et $\J_\infty^\cris$ et de 
schémas formels. Toutefois, on montre sans mal (en utilisant le même 
argument que dans la première partie de cette remarque) que les deux 
définitions coïncident.

\bigskip

Le résultat principal de \cite{breuil} est le suivant :

\begin{theo}
Le foncteur $\Mod$ réalise une anti-équivalence de catégories entre la 
catégorie des $\O_K$-groupes et la catégorie $\Mods$. De plus, cette 
anti-équivalence préserve les suites exactes courtes.
\end{theo}

\noindent
{\it Remarque.} Dans \emph{loc. cit.} Breuil donne une description
totalement explicite du quasi-inverse $\Gr$ du foncteur $\Mod$. Nous 
n'aurons pas besoin de cette description pour cet article et donc nous
ne la détaillons pas ici.

\section{Dualité sur les catégories de modules}

Cette section est consacrée à la définition de dualités (dans le sens 
\og anti-équivalence de catégories sur elle-même \fg) d'une part sur la 
catégorie $\Mods$ et d'autre part sur la catégorie des modules fortement 
divisibles.
On commence par traiter le cas des modules fortement divisibles pour 
lesquels la construction est simplifiée par l'existence de bases 
adaptées. On déduira ensuite la dualité sur la catégorie $\Mods$ de 
celle que l'on aura définie sur les modules fortement divisibles.

\subsection{Sur les modules fortement divisibles}
Soit $\calM$ un module fortement divisible. Le dual de $\calM$ est 
l'objet $(\calM^\vee, \Fil^1 \calM^\vee, \phi_1^\vee)$ défini de la 
façon suivante :
\begin{dinglist}{43}
\item $\calM^\vee = \hom_S (\calM, S)$ (où $\hom_S$ signifie que
l'on considère \emph{tous} les morphismes $S$-linéaires) ;
\item $\Fil^1 \calM^\vee$ est le sous-ensemble de $\calM^\vee$ formé
des $f : \calM \to S$ qui envoient $\Fil^1 \calM$ dans $\Fil^1 S$ ;
\item pour tout $f \in \Fil^1 \calM^\vee$, l'application $\phi_1^\vee
(f)$ est l'unique morphisme $S$-linéaire faisant commuter le
diagramme suivant :
$$\xymatrix @C=45pt {
\Fil^1 \calM \ar[r]^-{\phi_1} \ar[d]_f & \calM \ar[d]^{\phi_1^\vee (f)}
\\
\Fil^1 S \ar[r]^-{\phi_1} & S }$$
\end{dinglist}

\medskip

Pour que cette définition ait un sens, il faut montrer l'existence 
et l'unicité du morphisme $\phi_1^\vee (f)$. C'est l'objet du lemme 
suivant :

\begin{lemme}
\label{lem:transpphi}
Soient $\calM$ un module fortement divisible et $f : 
\Fil^1 \calM \to 
\Fil^1 S$ une application $S$-linéaire. Alors il existe une 
\emph{unique} application $g : \calM \to S$ faisant commuter le 
diagramme suivant :
$$\xymatrix @C=50pt {
\Fil^1 \calM \ar[r]^-{\Phi_1} \ar[d]_-{f} & \calM 
\ar[d]^-{g} \\
\Fil^1 S \ar[r]^-{\Phi_1} & S }$$
\end{lemme}

\begin{preuve}
L'unicité résulte simplement du fait que $\im \phi_1$ engendre $\calM$ 
en tant que $S$-module.

Pour l'existence, on considère $\pa{e_1, \ldots, e_d}$ est une base 
adaptée de $\calM$ pour les entiers $n_1, \ldots, n_d$ (voir lemme 
\ref{lem:adapt}). On pose $x_i = \phi_1(E(u)^{n_i} e_i) \in \calM$.
Du fait que $S$ est un anneau local, on vérifie que le lemme de Nakayama 
s'applique et implique que la famille des $x_i$ est génératrice. Comme 
en outre, elle a le bon cardinal, on en déduit qu'elle forme une 
$S$-base de $\calM$. Cette constatation permet de définir $g$ comme
l'unique application $S$-linéaire vérifiant $g(x_i) = \phi_1 \circ f 
(E(u)^{n_i} e_i)$ pour tout $i$. On vérifie alors facilement qu'elle
convient.
\end{preuve}

\bigskip

Nous voulons à présent montrer que le triplet $(\calM^\vee, \Fil^1 
\calM^\vee, \phi_1^\vee)$ définit un module fortement divisible. Il n'y 
a aucune difficulté à la vérification de la relation de compatibilité :
$$\phi^\vee_1\pa  {sx} = \frac 1 {c^r} \phi_1^\vee \pa{s} 
\phi_1^\vee\pa{E\pa u x}$$
pour tout $s \in \Fil^1 S$ et tout $x \in \calM^\vee$. Par 
ailleurs le fait que $\calM$ soit libre sur $S$ est évident. Les deux 
dernières propriétés \og $\calM^\vee / \Fil^1 \calM^\vee$ sans 
$p$-torsion \fg\ et \og $\phi_1(\Fil^1 \calM^\vee)$ engendre 
$\calM^\vee$ \fg\ résultent l'une comme l'autre du lemme suivant :

\begin{lemme}
Soit $\pa{e_1, \ldots, e_d}$ une base adaptée de $\calM$ pour 
les entiers $n_1, \ldots, n_d$. Alors la base duale $\pa{e^\vee_1, 
\ldots, e^\vee_d}$ de $\calM^\vee$ est également adaptée pour les entiers 
$n^\vee_1, \ldots, n^\vee_d$ avec $n^\vee_i = 1 - n_i$.

De plus si l'on pose $x_i = \Phi_1 \pa{E\pa u ^{n_i} e_i}$ et $x^\vee_i 
= \Phi_1^\vee (E\pa u ^{n^\vee_i} e^\vee_i)$, les familles $\pa{x_1, 
\ldots, x_d}$ et $\pa{x^\vee_1, \ldots, x^\vee_d}$ sont des bases duales 
l'une de l'autre.
\end{lemme}

\begin{preuve}
La famille des $e^\vee_j$ est définie par les égalités $e^\vee_j\pa{e_i} 
= \delta_{ij}$ où $\delta$ désigne le symbole de Kronecker. Soit $f \in 
\Fil^1 \calM^\vee$. On peut décomposer $f$ sur la base des $e^\vee_j$ et 
donc écrire :
$$f = s_1 e^\vee_1 + \ldots + s_d e^\vee_d$$
pour certains éléments $s_i \in S$. Comme $f \in \Fil^1 \calM^\vee$, 
il vient $f\pa{E\pa u ^{n_i} e_i} \in \Fil^1 S$. Or $f\pa{E\pa u 
^{n_i} e_i} = E\pa u^{n_i} s_i$. Cela prouve que $s_i$ est un élément de 
$\Fil^{n_i^\vee} S$ et donc que la base $\pa{e^\vee_1, \ldots, 
e^\vee_d}$ est adaptée pour les entiers $n^\vee_1, \ldots, n^\vee_d$.

\medskip

Passons à la seconde partie du lemme. On a déjà vu (dans la preuve 
du lemme \ref{lem:transpphi}) que $(x_1, \ldots, x_d)$ est une base de 
$\calM$. Considérons le diagramme commutatif suivant :
$$\xymatrix @C=50pt {
\Fil^1 \calM \ar[r]^-{\Phi_1} \ar[d]_-{f} & \calM 
\ar[d]^-{\Phi_r^\vee \pa f} \\
\Fil^1 S \ar[r]^-{\Phi_1} & S }$$
En prenant $f = E\pa u^{n_j^\vee} e^\vee_j$ et en regardant quelle est
l'image de $E\pa u^{n_i} e_i$ par chacun des deux chemins, on obtient
$x_j^\vee \pa{x_i} = \delta_{ij}$, ce qui conclut.
\end{preuve}

\bigskip

On déduit de ce qui précède la proposition suivante :

\begin{prop}
Le triplet $(\calM^\vee, \Fil^1 \calM^\vee, \phi_1^\vee)$ défini 
précédemment est un module fortement divisible.
\end{prop}

\bigskip

Si $f$ est un morphisme entre modules fortement divisibles, on vérifie 
sans peine que sa transposée (au sens classique) est compatible à toutes 
les structures et donc aussi un morphisme entre modules fortement 
divisibles. On a ainsi défini un foncteur ${}^\vee$ contravariant de la 
catégorie des modules fortement divisibles sur elle-même.

\begin{prop}
\label{prop:dualfortdiv}
Le foncteur ${}^\vee$ est une anti-équivalence de catégories 
(\emph{i.e.} une dualité). De plus, il transforme suites exactes courtes 
en suites exactes courtes.
\end{prop}

\begin{preuve}
Le foncteur ${}^\vee$ est son propre quasi-inverse. En effet, on a un 
morphisme canonique de $S$-modules de $\calM$ dans $\calM^{\vee\vee}$ 
donné par $x \mapsto \pa{f \mapsto f\pa x}$. C'est un isomorphisme 
puisque $\calM$ est libre sur $S$. Il est facile de vérifier que cet 
isomorphisme respecte $\Fil^1$ et $\phi_1$. Le seul point délicat est de 
montrer la surjectivité de $\Fil^1 \calM \to \Fil^1 \calM^{\vee\vee}$. 
Cela revient à montrer que si $x \not\in \Fil^1 \calM$, alors il existe 
$f \in \Fil^1 \calM^\vee$ tel que $f(x) \not\in \Fil^1 S$. On 
considère pour cela $(e_1, \ldots, e_d)$ une base adaptée de $\calM$ 
pour les entiers $n_1, \ldots, n_d$. On a alors une écriture :
$$x = s_1 e_1 + \cdots + s_d e_d$$
avec $s_i \in S$ et comme $x \not\in \Fil^1 \calM$, il existe un 
indice $i$ tel que $s_i \not\in \Fil^{n_i} S$. On vérifie alors aisément 
que l'application $S$-linéaire $f : \calM \to S$ définie par $f(e_i) = 
E(u)^{1-n_i}$ et $f(e_j) = 0$ pour $j \neq i$ convient.

\medskip

Montrons à présent le second point.
Considérons $0 \to \calM' \to \calM \to \calM'' \to 0$ une suite exacte
de modules fortement divisibles. Puisque $\calM''$ est un 
$S$-module libre, la suite $0 \to \calM''^\vee \to \calM^\vee \to 
\calM'^\vee \to 0$ est exacte comme suite de $S$-modules. Par 
ailleurs, la suite :
$$\xymatrix {
0 \ar[r] & \Fil^1 \calM''^\vee \ar[r] & \Fil^1 \calM^\vee \ar[r] & 
\Fil^1 \calM'^\vee \ar[r] & 0 }$$
est exacte à gauche. Il s'agit simplement de démontrer la surjectivité
de $\Fil^1 \calM^\vee \to \Fil^1 \calM'^\vee$. Notons $\Fil^1 \calN$ 
l'image de ce morphisme. On vérifie directement que ce $\Fil^1$ fait de 
$\calM'^\vee$ un second module fortement divisible : notons le $\calN$. 
On a une suite exacte $0 \to \calM''^\vee \to \calM^\vee \to \calN \to 
0$.

D'autre part, l'identité fournit un morphisme $f : \calN \to \calM'$ 
dans la catégorie $\pModS$ dont le dual s'insère dans le diagramme 
commutatif suivant :
$$\xymatrix {
0 \ar[r] & \Fil^1 \calM' \ar[r] \ar[d]_{f^\vee} & \Fil^1 \calM \ar[r] 
\ar@{=}[d] & \Fil^1 \calM'' \ar@{=}[d] \\
0 \ar[r] & \Fil^1 \calN^\vee \ar[r] & \Fil^1 \calM \ar[r] & \Fil^1 
\calM'' }$$
On en déduit que $f^\vee$ est un isomorphisme de modules fortement 
divisibles, puis qu'il en est de même de $f = (f^\vee)^\vee$. Finalement 
$\Fil^1 \calN = \Fil^1 \calM'$, et la surjectivité de $\Fil^1 \calM^\vee 
\to \Fil^1 \calM'^\vee$ en découle.
\end{preuve}

\subsection{Sur la catégorie $\Mods$}
\subsubsection{Définition de l'objet dual}

On introduit les deux $S$-modules $S_{K_0}$ et  $S_\infty$ définis par
$S_{K_0} = S \otimes_W K_0$ et $S_\infty = S \otimes_W K_0/W$. On pose 
$\Fil^1 S_{K_0} = \Fil^1 S \otimes_W K_0 \subset S_{K_0}$ (car $K_0$ est 
plat sur $W$). On a une projection $S_{K_0} \to S_\infty$ et on note 
$\Fil^1 S_\infty$ l'image de $\Fil^1 S_{K_0}$ dans $S_\infty$. On 
vérifie que $\phi_1 : \Fil^1 S \to S$ induit des applications $\Fil^1 
S_{K_0} \to S_{K_0}$ et $\Fil^1 S_\infty \to S_\infty$ que l'on
appelle encore $\phi_1$.

\medskip

Soit $\calM$ un objet de $\Mods$. Le dual de $\calM$ est l'objet 
$(\calM^\vee, \Fil^1 \calM^\vee, \phi_1^\vee)$ défini de la façon 
suivante :
\begin{dinglist}{43}
\item $\calM^\vee = \hom_S (\calM, S_\infty)$ (où $\hom_S$ signifie que
l'on considère \emph{tous} les morphismes $S$-linéaires) ;
\item $\Fil^1 \calM^\vee$ est le sous-ensemble de $\calM^\vee$ formé
des $f : \calM \to S_\infty$ qui envoient $\Fil^1 \calM$ dans $\Fil^1 
S_\infty$ ;
\item pour tout $f \in \Fil^1 \calM^\vee$, l'application $\phi_1^\vee
(f)$ est l'unique morphisme $S$-linéaire faisant commuter le
diagramme suivant :
$$\xymatrix @C=45pt {
\Fil^1 \calM \ar[r]^-{\phi_1} \ar[d]_f & \calM \ar[d]^{\phi_1^\vee (f)}
\\
\Fil^1 S_\infty \ar[r]^-{\phi_1} & S_\infty }$$
\end{dinglist}

\bigskip

Il reste à prouver l'existence et l'unicité du morphisme $\phi_1^\vee 
(f)$ mentionné ci-dessus. Comme dans le cas des modules fortement 
divisibles, l'unicité résulte simplement du fait que $\phi_1(\Fil^1 
\calM)$ engendre tout $\calM$. L'existence par contre est plus délicate. 
Nous consacrons tout le paragraphe suivant à son établissement.

\subsubsection{Existence du morphisme $\phi_1^\vee(f)$}

On commence par prouver deux lemmes :

\begin{lemme}
\label{lem:fortquot}
Pour tout objet $\calM$ de $\Mods$, il existe $\hat \calM$ et $\hat 
\calM'$ des modules fortement divisibles et une suite exacte :
$$\xymatrix {
0 \ar[r] & \hat \calM' \ar[r] & \hat \calM \ar[r] & \calM \ar[r] & 0 }$$
dans la catégorie $\pModS$.
\end{lemme}

\begin{preuve}
Par les alinéas 3.1.1 et 3.3.12 de \cite{bbm}, il existe des groupes 
$p$-divisibles $\H$ et $\H'$ et une suite exacte $0 \to \G \to \H \to 
\H' \to 0$. Le lemme s'en déduit par l'anti-équivalence de Breuil (en 
travaillant dans un premier temps modulo $p^n$ puis en passant à la 
limite).
\end{preuve}

\bigskip

\noindent
{\it Remarque.} On peut également démontrer le lemme précédent avec 
seulement des considérations d'algèbre linéaire.

\medskip

\begin{lemme}
\label{lem:fortquotdual}
Si $\calM$, $\hat \calM$ et $\hat \calM'$ sont comme dans le lemme 
\ref{lem:fortquot}, alors on a également une suite exacte :
$$\xymatrix {
0 \ar[r] & \hat \calM^\vee \ar[r] & \hat \calM'^\vee \ar[r] & \calM 
\ar[r] & 0 }$$
dans la catégorie $\pModS$.
\end{lemme}

\begin{preuve}
Définissons dans un premier temps les flèches. La première est 
simplement la transposée de l'inclusion $\hat \calM' \toinj \hat \calM$.
Pour la seconde, considérons $f \in \hat \calM'^\vee$ et $x \in \calM$. 
Soit $\hat x \in \hat \calM$ un relèvement de $x$. Comme $\calM$ est 
tué par une puissance de $p$, il existe un entier $n$ tel que $p^n \hat 
x \in \hat \calM'$. La réduction dans $S_\infty$ de $\frac 1 {p^n} 
f\pa{p^n \hat x}$ ne dépend ni de l'entier $n$, ni du relèvement $\hat 
x$ choisi. Cela permet de définir une application $S$-linéaire $\calM 
\to S_\infty$.

Par exactitude à gauche, le noyau de $\hat \calM^\vee \to \hat 
\calM'^\vee$ s'identifie à $\hom_S (\calM, S)$ qui est nul puisque 
$\calM$ est tué par une puissance de $p$. La première flèche est donc 
bien injective.

Prouvons l'exactitude au milieu. Soit $f : \calM' \to S$ une application 
$S$-linéaire. On suppose que l'image de $f$ dans $\calM^\vee$ est nulle 
et on veut montrer que $f$ se prolonge à $\calM$. Soient $x \in \calM$ 
et $n$ un entier tel que $p^n x \in \calM'$. Par hypothèse $\frac 1 
{p^n} f\pa{p^n x}$ est nul dans $S^\infty$, ce qui signifie que 
$f\pa{p^n x}$ est un multiple de $p^n$. On définit alors $f\pa x = 
\frac{f\pa{p^n x}}{p^n}$ (qui est bien défini puisque $S$ est intègre).

Passons à la surjectivité. Soit $f : \calM \to S_\infty$ une application 
$S$-linéaire. Considérons $\pa{\hat e_1, \ldots, \hat e_d}$ une base de 
$\hat \calM$ et notons $e_i$ l'image dans $\calM$ de $\hat e_i$. Pour 
tout $i$, notons $\hat x_i$ un relevé quelconque dans $S_{K_0}$ de $x_i 
= f\pa{e_i} \in S_\infty$. On définit  $\hat f \pa{\hat e_i} = \hat x_i$ 
et par linéarité on étend $\hat f$ en une application $\calM^\vee \to 
S_{K_0}$. On vérifie alors que la restriction de $\hat f$ à $\hat 
\calM'$ tombe dans $S$ et qu'elle induit $f$ dans $\calM^\vee$.

\medskip

Il reste à prouver l'exactitude au niveau des $\Fil^1$. L'injectivité et 
l'exactitude au milieu se traitent comme précédemment. Pour la 
surjectivité, considérons $f \in \Fil^1 \calM^\vee$ et $\pa{\hat e_1, 
\ldots, \hat e_d}$ une base adaptée de $\hat \calM$ pour les entiers 
$n_1, \ldots, n_d$. Notons $e_i$ l'image dans $\calM$ de $\hat e_i$. Par 
hypothèse $E(u)^{n_i} f(e_i) \in \Fil^1 S_\infty$ et cela assure qu'il 
existe $\hat x_i \in S_{K_0}$ relevant $f(e_i)$ et tel que $E(u)^{n_i} 
\hat x_i \in \Fil^1 S_{K_0}$. L'application $\hat f$ définie par $\hat 
f(\hat e_i) = \hat x_i$ se restreint alors à $\hat \calM'$ en un élément 
de $\Fil^r \hat \calM'^\vee$ qui est un antécédent de $f$.
\end{preuve}

\bigskip

On déduit finalement simplement des deux lemmes précédents l'existence 
de $\phi_1^\vee (f)$, pour $f \in \Fil^1 \calM^\vee$. En effet, par le 
lemme \ref{lem:fortquotdual}, $f$ se relève en un élément $\hat f \in
\Fil^1 \hat \calM'^\vee$ et on vérifie sans difficulté que l'image
de $\phi_1^\vee (\hat f)$ dans $\calM^\vee$ convient.

\subsubsection{Propriétés du foncteur de dualité}

Nous n'avons toujours pas montré que le triplet $(\calM^\vee, \Fil^1 
\calM^\vee, \phi_1^\vee)$ reste un objet de $\Mods$, mais pour cela
nous allons avoir besoin de l'exactitude que nous prouvons dans un 
premier temps. La preuve est basée sur une légère généralisation
du lemme \ref{lem:fortquot} que nous donnons ci-dessous :

\begin{lemme}
\label{lem:fortquotexact}
Soit $0 \to \calM \to \calX \to \calN \to 0$ une suite exacte dans la 
catégorie $\Mods$. Il existe des modules fortement divisibles $\hat 
\calM$, $\hat \calM'$, $\hat \calX$, $\hat \calX'$, $\hat \calN$ et 
$\hat \calN'$ qui s'insèrent dans le diagramme commutatif suivant :
$$\xymatrix @C=40pt @R=15pt {
& 0 \ar[d] & 0 \ar[d] & 0 \ar[d] \\
0 \ar[r] & \hat \calM' \ar[r] \ar[d] & \hat \calX' \ar[r] \ar[d] &
\hat \calN' \ar[r] \ar[d] & 0 \\
0 \ar[r] & \hat \calM \ar[r] \ar[d] & \hat \calX \ar[r] \ar[d] &
\hat \calN \ar[r] \ar[d] & 0 \\
0 \ar[r] & \calM \ar[r] \ar[d] & \calX \ar[r] \ar[d] & \calN \ar[r] 
\ar[d] & 0 \\
& 0 & 0 & 0 \\ }$$
où toutes les lignes et colonnes sont des suites exactes dans la 
catégorie $\pModS$.
\end{lemme}

\begin{preuve}
Par le lemme \ref{lem:fortquot}, on peut déjà construire un diagramme :
$$\xymatrix @C=40pt @R=15pt {
& 0 \ar[d] & & 0 \ar[d] \\
& \hat \calM' \ar[d] & & \hat \calN' \ar[d] \\
& \hat \calM \ar[d] & & \hat \calN \ar[d] \\
0 \ar[r] & \calM \ar[r] \ar[d] & \calX \ar[r] & \calN \ar[r] 
\ar[d] & 0 \\
& 0 & & 0 \\ }$$
Pour le reste, on pose $\hat \calX = \hat \calM \oplus \hat \calN$ et 
$\Fil^1 \hat \calX = \Fil^1 \hat \calM \oplus \Fil^1 \hat \calN$. Tout 
d'abord, on a $\hat \calX / \Fil^1 \hat \calX \simeq \hat \calM / \Fil^1 
\hat \calM \oplus \hat \calN / \Fil^1 \hat \calN$ ; c'est donc un module 
sans $p$-torsion. Considérons ensuite $(\hat e_1, \ldots, \hat e_d)$ une 
base adaptée de $\hat \calN$ pour les entiers $n_1, \ldots, n_d$. Notons 
$e_i$ l'image de $\hat e_i$
dans $\calN$ et $e'_i \in \calX$ un relevé de $e_i$. Si $n_i = 0$, on a 
$e_i \in \Fil^1 \calN$ et donc on peut choisir $e'_i \in \Fil^1 \calX$, 
ce que l'on ne se prive pas de faire. On définit une application 
surjective $f : \hat \calX \to \calX$ qui coïncide sur $\hat \calM$ avec 
la projection $\hat \calM \to \calM$ et qui est définie sur $\hat \calN$ 
par les égalités $f(0 \oplus \hat e_i) = e'_i$. On vérifie que $f$  
induit une application surjective $\Fil^1 \hat \calX \to \Fil^1 \calX$
qui fait commuter le diagramme suivant :
$$\xymatrix @C=40pt @R=15pt {
0 \ar[r] & \hat \calM \ar[r] \ar@{->>}[d] & \hat \calX \ar[r] 
\ar@{->>}[d] & \hat \calN \ar[r] \ar@{->>}[d] & 0 \\
0 \ar[r] & \calM \ar[r] & \calX \ar[r] & \calN \ar[r] & 0 }$$

On définit à présent l'opérateur $\phi_1 : \Fil^1 \hat \calX \to \calX$ 
de la façon suivante. Sur $\Fil^1 \hat \calM$, il coïncide avec 
l'application
$\phi_1 : \Fil^1 \hat \calM \to \calM$. Il ne reste qu'à donner ses
valeurs sur les éléments $E(u)^{n_i} (0 \oplus \hat e_i)$. 
Notons $x_i = \phi_1 (E(u)^{n_i}e_i) \in \calN$ et $x'_i \in \calX$ un 
relevé un $x_i$. Les éléments $x'_i$ et $\phi_1 (E(u)^{n_i} \hat e_i)$ 
s'envoient tous deux sur $x_i$ dans $\calN$ ; ils admettent donc un 
antécédent commun dans $\hat \calX$, disons $\hat x_i$. On pose
$\phi_1(E(u)^{n_i} (0 \oplus \hat e_i)) = \hat x_i$ et on vérifie que 
l'on obtient bien ainsi un module fortement divisible. De plus, 
l'application $f$ ainsi que la projection canonique $\hat \calX \to \hat 
\calN$ sont par construction compatibles à $\phi_1$.

On pose $\hat \calX' = \ker f$ et $\Fil^1 \hat \calX' = 
\hat \calX' \cap \Fil^1 \hat \calX$. Le quotient $\hat \calX' / \Fil^1 
\hat \calX'$ s'injecte dans $\hat \calX / \Fil^1 \hat \calX$ et est donc 
également sans $p$-torsion. Par ailleurs, l'application $\phi_1$ 
précédemment définie induit une flèche $\Fil^1 \hat \calX' \to \calX'$ 
qui s'insère dans le diagramme commutatif suivant :
$$\xymatrix @C=40pt @R=15pt {
0 \ar[r] & \Fil^1 \hat \calM' \ar[r] \ar[d]_-{\phi_1} & \Fil^1 \hat 
\calX' \ar[r] \ar[d]_-{\phi_1} & \Fil^1 \hat \calN' \ar[r] 
\ar[d]_-{\phi_1} & 0 \\
0 \ar[r] & \hat \calM' \ar[r] & \hat \calX' \ar[r] & \hat \calN'
\ar[r] & 0 
}$$
où les deux lignes sont exactes par application du lemme du serpent. On 
en déduit que $\phi_1(\Fil^1 \hat \calX')$ engendre $\hat \calX'$ et que 
$\hat \calX'$ est un $S$-module libre. Il s'agit donc d'un module 
fortement divisible et cela termine la démonstration.
\end{preuve}

\begin{lemme}
\label{lem:exact}
Le foncteur ${}^\vee$ conserve les suites exactes courtes.
\end{lemme}

\begin{preuve}
On ne traite que l'exactitude en tant que $S$-module, celle au niveau 
des $\Fil^1$ étant en tout point analogue. Soit $0 \to \calM \to \calX 
\to \calN \to 0$ une suite exacte dans $\Mods$. D'après le lemme
\ref{lem:fortquotexact}, il existe un diagramme commutatif de la forme :
$$\xymatrix @C=40pt @R=15pt {
& 0 \ar[d] & 0 \ar[d] & 0 \ar[d] \\
0 \ar[r] & \hat \calM' \ar[r] \ar[d] & \hat \calX' \ar[r] \ar[d] &
\hat \calN' \ar[r] \ar[d] & 0 \\
0 \ar[r] & \hat \calM \ar[r] \ar[d] & \hat \calX \ar[r] \ar[d] &
\hat \calN \ar[r] \ar[d] & 0 \\
0 \ar[r] & \calM \ar[r] \ar[d] & \calX \ar[r] \ar[d] & \calN \ar[r] 
\ar[d] & 0 \\
& 0 & 0 & 0 \\ }$$
où $\hat \calM$, $\hat \calM'$, $\hat \calX$, $\hat \calX'$, 
$\hat \calN$ et $\hat \calN'$ sont des modules fortement divisibles. En 
dualisant, on obtient :
$$\xymatrix @C=40pt @R=15pt {
0 \ar[r] & \hat \calN'^\vee \ar[r] \ar@{->>}[d] & \hat \calX'^\vee 
\ar[r] \ar@{->>}[d] & \hat \calM'^\vee \ar[r] \ar@{->>}[d] & 0 \\
0 \ar[r] & \calN^\vee \ar[r] & \calX^\vee \ar[r] & \calM ^\vee } $$
où la ligne du haut est exacte (proposition \ref{prop:dualfortdiv}) 
et les flèches verticales sont surjectives (lemme 
\ref{lem:fortquotdual}). La surjectivité $\calX^\vee \to \calM^\vee$ 
résulte alors d'une chasse au diagramme triviale et permet de conclure.
\end{preuve}

\medskip

La propriété suivante résume finalement les propriétés du foncteur 
${}^\vee$ :

\begin{prop}
\label{5:prop:dualmr}
Le foncteur ${}^\vee$ est une dualité de la catégorie $\Mods$. De plus, 
il transforme suites exactes courtes en suites exactes courtes.
\end{prop}

\begin{preuve}
Il ne reste qu'à prouver que le dual d'un objet $\calM$ de $\Mods$ est 
encore un objet de $\Mods$. Par le lemme \ref{lem:exact}, il suffit de 
le faire lorsque $\calM$ est tué par $p$.

Dans ce cas, $\calM$ est un $S/pS$-module libre et de $\hom_S (S/pS, 
S_\infty) = S/pS$, on déduit que $\calM^\vee$ est aussi un $S/pS$-module 
libre. Il ne reste qu'à prouver que $\phi_1(\Fil^1 \calM^\vee)$ engendre 
$\calM^\vee$. Mais si $0 \to \hat \calM' \to \hat \calM \to \calM \to 0$ 
est une suite exacte dans $\pModS$ avec $\hat \calM$ et $\hat \calM'$ 
des modules fortement divisibles, on a le diagramme commutatif suivant :
$$\xymatrix @C=50pt {
\Fil^1 \hat \calM'^\vee \ar@{->>}[r] \ar[d]_-{\phi_1} & 
\Fil^1 \calM^\vee \ar[d]^-{\phi_1} \\
\hat \calM'^\vee \ar@{->>}[r] & \calM^\vee }$$
et $\phi_1(\Fil^1 \hat \calM'^\vee)$ engendre $\hat \calM'^\vee$ puisque 
$\hat \calM'^\vee$ est encore un module fortement divisible. Une chasse 
au diagramme permet alors de conclure.
\end{preuve}

\section{Dualité de Cartier}

Dans cette dernière section, on montre que la dualité définie 
précédemment correspond \emph{via} le foncteur $\Mod$ à la dualité de
Cartier sur les schémas en groupes.

\subsection{Rappels de théorie de Dieudonné cristalline}
On considère $\G$ un $\O_K$-groupe. Comme en \ref{subsec:defmod}, on 
note $\G_n$ la réduction modulo $p^n$ de $\G$. C'est un schéma en 
groupes commutatifs fini et plat sur la base $T_n = \spec (\O_K/p^n)$. 
On rappelle par ailleurs que $E_n$ désigne le schéma $\spec (S_n)$ où 
$S_n = S/p^nS$. La projection $S_n \to \O_K/p^n$ qui envoie $u$ sur 
$\pi$ définit un épaississement $T_n \toinj E_n$.

\medskip

Avant de poursuivre, notons également que, par la suite, nous aurons à 
considérer les objets $\hom(\calF, \calF')$ et $\hhom(\calF, \calF')$
pour $\calF$ et $\calF'$ deux faisceaux de groupes abéliens sur un 
certain site. Le premier désignera toujours l'ensemble des morphismes
entre $\calF$ et $\calF'$ dans la catégorie des faisceaux de 
\emph{groupes abéliens}. Le second, quant à lui, est une version 
faisceautique du premier : si $U$ est un objet du site, on a par
définition :
$$\hhom(\calF, \calF')(U) = \hom(\calF_{|U}, \calF'_{|U})$$
où les morphismes sont toujours considérés dans la catégorie des
faisceaux abéliens.

De même, on définit les objets $\ext^1(\calF, \calF')$ et 
$\eext^1(\calF, \calF')$ qui sont respectivement un groupe abélien et un 
faisceau de groupes abéliens. On prêtera attention au fait qu'ici 
l'association :
$$U \mapsto \ext^1(\calF_{|U}, \calF'_{|U})$$
ne définit pas en général un faisceau mais seulement un préfaisceau. 
Par définition, $\eext^1(\calF, \calF')$ est le faisceau associé à ce 
préfaisceau.

\bigskip

Dans \cite{bbm} (voir définition 3.1.5), Berthelot, Breen 
et Messing associent à $\G_n$ un cristal sur le site\footnote{En 
réalité dans \cite{bbm}, il n'est pas du tout question de topologie 
syntomique. Cependant d'après le corollaire 2.3.11 de \cite{bbm}, il 
s'agit bien du même préfaisceau.} $(T_n/E_n)_\SYNCRIS$, appelé 
\emph{cristal de Dieudonné} de $\G_n$ et noté $\D(\G_n)$. Par 
définition, on a $\D(\G_n) = \eext^1 (\G_n, \O_{T_n/E_n})$.

\medskip

D'autre part, de façon très générale, si $A$ et $B$ sont deux objets 
d'une catégorie abélienne tués par un entier $N$, on a une flèche 
canonique :
$$\ext^1 (A,B) \to \hom(A,B)$$
qui à une extension $E$ associe la flèche du serpent associée au 
diagramme commutatif à lignes exactes suivant :
$$\xymatrix @C=30pt @R=15pt {
0 \ar[r] & B \ar[r] \ar[d]^0 & E \ar[r] \ar[d]^N & A \ar[r] \ar[d]^0 & 0 
\\
0 \ar[r] & B \ar[r] & E \ar[r] & A \ar[r] & 0 }$$
Puisque $\G_n$ et $\O_{T_n/E_n}$ sont tués par $p^n$, ceci s'applique à 
notre situation et fournit une flèche canonique :
\begin{equation}
\label{5:eq:snake}
\sigma : \D(\G_n) \to \hhom(\G_n, \O_{T_n/E_n})
\end{equation}
qui, d'après la proposition 4.2.9 de \cite{bbm}, induit un isomorphisme 
sur les sections globales. Autrement dit :
\begin{equation}
\label{5:eq:breuilbbm}
\D(\G_n)(T_n, E_n) \simeq \hom(\G_n, \O_{T_n/E_n}) = \Mod(\G)
\end{equation}
la dernière égalité étant obtenue grâce à l'adjonction des foncteurs 
$w^\star$ et $w_\star$ (les $\hom$ calculés sur les petit et gros sites
syntomiques sont les mêmes puisque $\G_n$ est représentable par un
schéma syntomique sur $T_n$). L'isomorphisme (\ref{5:eq:breuilbbm}) est 
celui qui fournit le lien entre le point de vue de \cite{bbm} (cristal 
de Dieudonné) et le point de vue de \cite{breuil} (objet de $\Mods$).

\subsection{Construction du morphisme de comparaison}
\label{5:sec:morcomp}
On garde les notations du paragraphe précédent. On note, en outre, 
$\G^\vee$ le dual de Cartier de $\G$ et $\G_m^\vee = \G^\vee 
\times_{\spec(\O_K)} T_m$ pour $m=n$ et $m=n+1$. On désigne encore par
$\G_n^\vee$ le faisceau sur le site $(T_n/E_n)_\SYNCRIS$ défini par 
le schéma $\G_n^\vee$. On rappelle que $\G_n^\vee = \hhom(\G_n, 
\Gm)$ où $\Gm$ est défini par $\Gm(U,T) = \Gamma(U, \O_U)^\star$.

\medskip

Le but, ici, est d'obtenir un isomorphisme canonique et fonctoriel :
$$\Mod(\G)^\vee \to \Mod(\G^\vee).$$
Or, on dispose d'une suite exacte de faisceaux abéliens :
$$\xymatrix {
0 \ar[r] & 1 + \J_{T_n/E_n} \ar[r] & \O_{T_n/E_n}^\star \ar[r] & \Gm 
\ar[r] & 0 }$$
et la flèche de cobord associée au foncteur $\hom(\G_n, \cdot)$ induit 
un morphisme $\G_n^\vee \to \eext^1 (\G_n, 1 + \J_{T_n/E_n})$. Par 
ailleurs, on a un morphisme $\log : 1 + \J_{T_n/E_n} \to \O_{T_n/E_n}$ 
défini par :
$$\log(1+x) = 1 - x + \frac{x^2}2 - \frac{x^3}3 + \frac{x^4}4 - \cdots$$
qui induit, par fonctorialité, un morphisme $\eext^1 (\G_n, 1 + 
\J_{T_n/E_n}) \to \eext^1 (\G_n, \O_{T_n/E_n}) = \D(\G_n)$. En composant 
les deux flèches précédentes, on définit :
$$\alpha_\CRIS : \G_n^\vee \to \D(\G_n).$$
Finalement, par application du foncteur $\hom(\cdot, \O_{T_n/E_n})$, on 
obtient :
$$\alpha^\star_\CRIS : \hom_{\O_{T_n/E_n}} (\D(\G_n), \O_{T_n/E_n}) \to 
\hom (\G_n^\vee, \O_{T_n/E_n})$$
où la notation \og $\hom_{\O_{T_n/E_n}}$ \fg\ signifie que l'on se 
restreint aux morphismes $\O_{T_n/E_n}$-linéaires. D'après le théorème
5.2.7 de \cite{bbm}, $\alpha^\star_\CRIS$ est un isomorphisme.

Le but de $\alpha^\star_\CRIS$ s'identifie (grâce à l'adjonction des
foncteurs $w^\star$ et $w_\star$) à $\Mod(\G^\vee)$, tandis que la 
source est naturellement munie d'un morphisme $\gamma$ (obtenu en 
regardant les sections globales) vers $\hom_{S_n} (\Mod(\G), S_n)$.

\begin{lemme}
\label{5:lem:gamma}
Le morphisme $\gamma : \hom_{\O_{T_n/E_n}} (\D(\G_n), \O_{T_n/E_n}) \to 
\hom_{S_n} (\Mod(\G), S_n)$ est un isomorphisme.
\end{lemme}

\begin{preuve}
On remarque que $\D(\G_n)$ et $\O_{T_n/E_n}$ sont tous les deux des 
cristaux sur le site $(T_n/E_n)_\SYNCRIS$. Le lemme résulte alors de la 
description de la catégorie de ces cristaux en terme de modules à 
connexion intégrable et quasi-nilpotente (dans cette situation, la 
connexion est nécessairement nulle).
\end{preuve}

\bigskip

Au final, la composée $\beta^\star = \alpha^\star_\CRIS \circ 
\gamma^{-1}$ fournit un isomorphisme :
$$\beta^\star : \Mod(\G)^\vee \to \Mod(\G^\vee)$$
dont on vérifie directement qu'il est fonctoriel en $\G$. Il reste à 
prouver que $\beta^\star$ est isomorphisme \emph{dans la catégorie 
$\Mods$}. C'est l'objet des paragraphes suivants.

\subsection{Cas des groupes de la forme $\H(n)$}
Dans ce paragraphe, on se place dans le cas particulier où $\G$ est le 
noyau de la multiplication par $p^n$ sur un groupe $p$-divisible $\H$. 
Notons, pour simplifier, $\H_n = \H \times_{\spec(\O_K)} T_n$ et pour 
tout $m$, $\H_n(m)$ le noyau de la multiplication par $p^m$ sur $\H_n$. 
La supposition que l'on vient de faire entraîne alors $\G_n = \H_n(n)$.

\bigskip

Cette hypothèse supplémentaire a l'avantage de fournir un inverse au 
morphisme $\sigma$ défini en (\ref{5:eq:snake}). En effet, on dispose 
d'une suite exacte :
$$\xymatrix { 
0 \ar[r] & \G_n \ar[r] & \H_n(2n) \ar[r]^-{p^n} & \G_n \ar[r] & 0 }$$
et on note $s$ le morphisme de cobord :
$$s : \hhom(\G_n, \O_{T_n/E_n}) \to \eext^1 (\G_n, \O_{T_n/E_n}) = 
\D(\G_n)$$ 
associé au foncteur $\hom(\cdot, \O_{T_n/E_n})$.

\begin{lemme}
Les morphismes $\sigma$ et $s$ sont inverses l'un de l'autre. En 
particulier, le faisceau $\hhom(\G_n, \O_{T_n/E_n})$ est un cristal
sur $(T_n/E_n)_\SYNCRIS$.
\end{lemme}

\begin{preuve}
On vérifie facilement, en déroulant les définitions, que $\sigma 
\circ s = \id$. Il suffit, pour conclure, de prouver que $s$ est un
épimorphisme. Considérons $(U,T) \in (T_n/E_n)_\SYNCRIS$ avec $T = \spec 
S$ affine et le diagramme commutatif suivant :
$$\xymatrix @C=50pt {
\hhom(\G_n, \O_{T_n/E_n}) (T_n,E_n) \otimes_{S_n} S \ar[r] 
\ar[d]_{s_{(T_n,E_n)} \otimes \id} & \hhom(\G_n, \O_{T_n/E_n}) 
(U,T) \ar[d]_{s_{(U,T)}} \\
\D(\G_n) (T_n,E_n) \otimes_{S_n} S \ar[r] & \D(\G_n) (U,T) } $$
La flèche de gauche est un isomorphisme, ainsi que celle du bas puisque 
$\D(\G_n)$ est un cristal sur $(T_n/E_n)_\SYNCRIS$. On en déduit que
$s_{(U,T)}$ est surjectif, ce qui suffit pour conclure.
\end{preuve}

\bigskip

Notons :
$$\sigma^\star : \hom_{\O_{T_n/E_n}} (\hhom(\G_n, \O_{T_n/E_n}), 
\O_{T_n/E_n}) \to \hom_{\O_{T_n/E_n}} (\D(\G_n), \O_{T_n/E_n})$$
$$\text{(resp. }s^\star :  \hom_{\O_{T_n/E_n}} (\D(\G_n), \O_{T_n/E_n}) 
\to \hom_{\O_{T_n/E_n}} (\hhom(\G_n, \O_{T_n/E_n}), \O_{T_n/E_n})
\text{)}$$
le morphisme induit par $\sigma$ (resp. par $s$) \emph{via} le foncteur 
$\hom_{\O_{T_n/E_n}} (\cdot, \O_{T_n/E_n})$.

\subsubsection{Les morphismes $\beta_\CRIS$ et $\beta_\syn$}
\label{5:subsec:betasyn}

On pose $\beta_\CRIS = \sigma \circ \alpha_\CRIS$. Sur le 
\emph{petit} site syntomique, à partir de la suite exacte (que l'on 
déduit de (\ref{5:eq:exactsyn})) :
$$\xymatrix {
0 \ar[r] & 1 + \J_n^\cris \ar[r] & {\O_n^\cris}^\star \ar[r] & 
\O_n^\star \ar[r] & 0 } $$
on définit un morphisme $\hhom (\G_n, \O_n^\star) \to \eext^1 (\G_n, 1 + 
\J_n^\cris)$ qui fournit, après composition par $\log$, 
un morphisme $\alpha_\syn : \hhom (\G_n, \O_n^\star) \to \eext^1 (\G_n, 
\O_n^\cris)$. Le morphisme $\beta_\syn$ s'obtient de façon analogue en 
composant par la flèche canonique $\eext^1 (\G_n, \O_n^\cris) \to 
\hhom (\G_n, \O_n^\cris)$.

\medskip

On vérifie, en déroulant les définitions, que $\beta_\syn$ se décrit 
localement de la façon explicite suivante. Soit $f \in \hhom (\G_n, 
\O_n^\star)$. Soit $U \in (T_n)_\syn$. On cherche à décrire 
l'élément $\beta_\syn(f(U)) \in \hom ({\G_n}_{|U}, {\O_n^\cris}_{|U})$. 
Considérons pour cela $V \in U_\syn$ suffisamment petit pour que la 
suite $0 \to 1 + \J_n^\cris (V) \to {\O_n^\cris}^\star (V) \to
\O_n^\star (V) \to 0$
soit exacte. L'élément $f$ fournit par restriction à $V$ une application 
$g : \G_n(V) \to \O_n^\star (V)$. Soient $x \in \G_n(V)$ et $y \in 
{\O_n^\cris}^\star (V)$ un relevé quelconque de $g(x)$. On vérifie 
directement que $y^{p^n}$ ne dépend que de $x$ et que c'est un 
élément de $1 + \J_n^\cris (V)$. Le morphisme $\beta_\syn(f(U))$ est 
alors celui qui, sur $V$, associe $\log(y^{p^n})$ à $x$.

\medskip

Finalement, notons que l'on aura également besoin d'utiliser le 
morphisme $i_\star \beta_\syn$ où on rappelle que $i$ désigne 
l'inclusion $T_n \toinj T_{n+1}$ ; par abus, on notera ce morphisme
encore $\beta_\syn$.

\bigskip

Si $\calF$ est un faisceau de $(T_n)_\SYN$ et $\calF'$ un 
faisceau de $(T_n/E_n)_\SYNCRIS$, les morphismes d'adjonction 
permettent de construire un morphisme canonique :
\begin{equation}
\label{5:eq:adjfais}
w^\star \hhom (\calF, w_\star \calF') \to \hhom (w^\star \calF, 
\calF').
\end{equation}
Avec $\calF = \G_n$ et $\calF' = \O_{T_n/E_n}$, on obtient une flèche
$w^\star \hhom (\G_n, \O_n^\cris) \to \hhom (\G_n, \O_{T_n/E_n})$
puis, par application du foncteur $\hom(\cdot, \O_{T_n/E_n})$ :
$$\begin{array}{lcl}
\hom_{\O_{T_n/E_n}} ( \hhom (\G_n, \O_{T_n/E_n}), 
\O_{T_n/E_n}) & \to & \hom_{\O_{T_n/E_n}}( w^\star \hhom (\G_n, 
\O_n^\cris), \O_{T_n/E_n}) \\ \vphantom{\bigg(} & &
\hspace{2em} \simeq \hom_{\O_n^\cris} (\hhom (\G_n, \O_n^\cris), 
\O_n^\cris)
\end{array} $$
et finalement par restriction au petit site puis application du foncteur 
$i_\star$, un morphisme :
$$\gamma_1 : \hom_{\O_{T_n/E_n}} ( \hhom (\G_n, \O_{T_n/E_n}),
\O_{T_n/E_n}) \to \hom_{\O_n^\cris} (\hhom (\G_n, \O_n^\cris),
\O_n^\cris)$$
où, cette fois-ci, le dernier $\hom$ est calculé sur le 
$(T_{n+1})_\syn$. (On remarque que puisque $\G_n$ est représentable par 
un schéma syntomique, on a l'identification $i_\star \hhom 
(\G_n, \O_n^\cris) = \hhom (\G_n, i_\star \O_n^\cris)$.)

\medskip

De plus, on vérifie facilement que si $\calF'$ est de la forme $w^\star 
\calF''$, le morphisme (\ref{5:eq:adjfais}) est un isomorphisme. 
Autrement dit $\hhom (w^\star \calF, w^\star \calF') = w^\star \hhom 
(\calF, \calF')$ pour $\calF$ et $\calF'$ des faisceaux sur le gros site 
syntomique. En particulier, en prenant $\calF = \G_n$ et $\calF' = 
\O_n^\star$, on obtient un isomorphisme (sur le site cristallin) entre 
$\G_n^\vee$ et $w^\star \hhom(\G_n, 
\O_n^\star)$. On en déduit que $\hhom(\G_n, \O_n^\star)$ est le faisceau 
sur $(T_n)_\SYN$ défini par le schéma $\G_n^\vee$. On note ce faisceau 
encore $\G_n^\vee$. Les propriétés d'adjonction fournissent des 
égalités :
$$\hom (\G_n^\vee, \O_{T_n/E_n}) = \hom (\G_n^\vee, \O_n^\cris) = \Mod( 
\G^\vee)$$
les $\hom$ étant calculés sur les gros ou petits sites.

\medskip

Le diagramme commutatif suivant résume les liens entre nombreux des 
morphismes introduits jusqu'alors :
$$\xymatrix @C=60pt {
\hom_{\O_{T_n/E_n}} (\D(\G_n), \O_{T_n/E_n}) \ar@<4pt>[d]_-{\sigma^\star 
\hspace{7pt}} \ar[rd]^-{\alpha_\CRIS^\star}_-{\sim} 
\ar@/_4.5cm/[ddd]_{\sim}^{\gamma} \\
\hom_{\O_{T_n/E_n}} (\hhom(\G_n, \O_{T_n/E_n}), \O_{T_n/E_n}) 
\ar@<4pt>[u]_-{\hspace{7pt} s^\star} 
\ar[r]^-{\beta_\CRIS^\star}_-{\sim} \ar[d]^-{\gamma_1} & \hom 
(\G_n^\vee, \O_{T_n/E_n}) \ar@{=}[d] \\ 
\hom_{\O_n^\cris} (\hhom(\G_n, \O_n^\cris), \O_n^\cris) 
\ar[r]^-{\beta_\syn^\star} \ar[d]^-{\gamma_2} & \hom (\G_n^\vee, 
\O_n^\cris) \ar@{=}[d] \\ 
\Mod(\G)^\vee \ar[r]^-{\beta^\star}_-{\sim} & \Mod(\G^\vee) }$$
Le morphisme $\gamma_2$ est obtenu simplement en regardant le morphisme 
induit sur les sections globales et les morphismes $\beta_\CRIS^\star$ 
et $\beta_\syn^\star$ ont des définitions évidentes. Notons de plus 
que tous les faisceaux syntomiques sont considérés sur le site 
$(T_{n+1})_\syn$. Comme $\beta^\star_\syn \circ \gamma_1$ est un 
isomorphisme, le morphisme $\gamma_1$ est injectif. Notons 
$\hom^\CRIS_{\O_n^\cris} (\hhom(\G_n, \O_n^\cris), \O_n^\cris)$ son 
image. Le diagramme précédent se modifie alors de la façon suivante :
$$\xymatrix @C=60pt {
\hom_{\O_{T_n/E_n}} (\D(\G_n), \O_{T_n/E_n}) \ar@<4pt>[d]_-{\sigma^\star 
\hspace{7pt}} \ar[rd]^-{\alpha_\CRIS^\star}
\ar@/_4.5cm/[ddd]_{\gamma} \\
\hom_{\O_{T_n/E_n}} (\hhom(\G_n, \O_{T_n/E_n}), \O_{T_n/E_n}) 
\ar@<4pt>[u]_-{\hspace{7pt} s^\star} 
\ar[r]^-{\beta_\CRIS^\star} \ar[d]_-{\gamma_1} & \hom 
(\G_n^\vee, \O_{T_n/E_n}) \ar@{=}[d] \\ 
\hom^\CRIS_{\O_n^\cris} (\hhom(\G_n, \O_n^\cris), \O_n^\cris) 
\ar[r]^-{\beta_\syn^\star} \ar[d]_-{\gamma_2} & \hom (\G_n^\vee, 
\O_n^\cris) \ar@{=}[d] \\ 
\Mod(\G)^\vee \ar[r]^-{\beta^\star} & \Mod(\G^\vee) }$$
où désormais \emph{toutes} les flèches sont des isomorphismes comme on 
le vérifie facilement. Lorsque nous aurons à considérer par la suite des
morphismes $\gamma_1$, $\gamma_2$ et $\beta^\star_\syn$, l'ensemble de 
départ ou d'arrivée (selon le cas) sera toujours 
$\hom^\CRIS_{\O_n^\cris} (\hhom(\G_n, \O_n^\cris), \O_n^\cris)$. En 
particulier, ces trois morphismes deviennent des isomorphismes.

\subsubsection{Compatibilité à $\Fil^1$}

Le but de ce paragraphe est de montrer que le morphisme $\beta^\star$ 
envoie $\Fil^1 \Mod(\G)^\vee$ sur $\Fil^1 \Mod(\G^\vee)$. En réalité, 
cela résulte presque directement du lemme suivant :

\begin{lemme}
\label{5:lem:gammafil}
Soit $f \in \Fil^1 \Mod(\G)^\vee$. La restriction de $\gamma_2^{-1}(f)$ 
à $\hhom (\G_n, \J_n^\cris)$ tombe dans $\J_n^\cris$.
\end{lemme}

\begin{preuve}
Pour cette preuve on travaille sur le \emph{petit} site 
cristallin $(T_n/E_n)_\syncris$ : c'est la restriction du 
gros site aux couples $(U,T)$ pour lesquels $U$ est syntomique sur 
$T_n$. Notons $\tilde f$ la restriction du faisceau $\gamma_1^{-1} \circ 
\gamma_2 ^{-1} (f)$ au petit site $(T_n/E_n)_\syncris$. Il suffit de 
montrer que la restriction de $\tilde f$ à $\hhom 
(\G_n, \J_{T_n/E_n})$ tombe dans $\J_{T_n/E_n}$ et pour cela de
construire un morphisme $\tilde g : \hhom (\G_n, \Ga) \to \Ga$ faisant 
commuter le diagramme suivant :
\begin{equation}
\label{5:eq:tildeg}
\raisebox{0.5\depth} { \xymatrix {
0 \ar[r] & \hhom (\G_n, \J_{T_n/E_n}) \ar[r] & \hhom (\G_n, 
\O_{T_n/E_n}) \ar[r]^-{\pr_\star} \ar[d]_{\tilde f} & \hhom(\G_n, \Ga) 
\ar@{.>}[d]^-{\tilde g} \\
0 \ar[r] & \J_{T_n/E_n} \ar[r] & \O_{T_n/E_n} \ar[r]^-{\pr} & \Ga 
\ar[r] & 0 } }
\end{equation}
Comme $\hhom(\G_n, \Ga) = w^\star \hhom(\G_n, \O_n)$ et $\Ga = w^\star 
\O_n$, il revient au même de construire un morphisme $\hhom (\G_n, \O_n) 
\to \O_n$ sur le petit site syntomique $(T_n)_\syn$. On construit 
celui-ci localement.

\medskip

D'après la proposition 3.2.9 de \cite{breuil} (en remplaçant $\H(1)$ par 
$\H(n)$ ce qui ne modifie pas la preuve) le morphisme $\hom(\G_n, 
\O_{T_n/E_n}) \to \hom(\G_n, \O_n)$ est surjectif, et donc l'hypothèse 
du lemme fournit, à partir de $f$, une flèche $g : \hom (\G_n, \O_n) 
\to \O_K/p^n$.

Par ailleurs, puisque $\G_n$ est fini sur $T_n$, il est nécessairement 
affine : notons $\calA_n$ son anneau et $c : \calA_n \to \calA_n 
\otimes_{O_K/p^n} \calA_n$ la comultiplication. Notons $\calA_n^c$ le 
noyau de $c - 1\otimes \id - \id \otimes 1$. Considérons $U \in 
(T_n)_\syn$, c'est en particulier un schéma plat sur $T_n$. On 
suppose en outre que $U = \spec R$ est un schéma affine. On a alors :
$$\hhom(\G_n, \O_n)(U) = \acco{ x \in \calA \otimes_{O_K/p^n} R \, / \, 
c(x) = 1 \otimes x + x\otimes 1 } = \calA_n^c \otimes_{O_K/p^n} R$$
la dernière égalité provenant de la platitude de $R$ sur $\O_K/p^n$.
Le morphisme $g : \calA_n^c \to \O_K/p^n$ construit précédemment donne,
par tensorisation par $R$ au dessus de $\O_K/p^n$, une flèche 
$\hhom(\G_n, \O_n)(U) \to \O_n(U)$. Par recollement, on construit un 
morphisme de faisceaux sur le site syntomique qui correspond à un 
morphisme $\tilde g : \hhom(\G_n, \Ga) \to \Ga$ sur le site 
cristallin.
Il ne reste plus qu'à vérifier que $\tilde g$ fait commuter le diagramme 
(\ref{5:eq:tildeg}). On considère pour cela $(U,T) \in 
(T_n/E_n)_\syncris$ et on vérifie la commutativité du diagramme sur 
l'ouvert $(U,T)$. On peut supposer $U$ et $T$ affines, disons $U = \spec 
R$ et $T = \spec S$. On considère le diagramme suivant :
$$\xymatrix @C=70pt {
\hhom (\G_n, \O_{T_n/E_n}) (U,T) \ar@{=}[d] \ar[r]^-{\pr_{\star (U,T)}} 
& \hhom (\G_n, \Ga) (U,T) \ar@{=}[d] \\
\Mod(\G_n) \otimes_{S_n} S \ar[d]_-{f \otimes \id} 
\ar[r]^-{\pr_{\star(T_n,E_n)} \otimes \pr_{(U,T)}} & \calA_n^c 
\otimes_{O_K/p^n} R \ar[d]^-{g \otimes \id} \\
S \ar[r]^-{\pr_{(U,T)}} & R }$$
On vérifie directement que le carré du haut est commutatif. Celui du bas 
l'est également par construction. Ainsi tout le diagramme commute, ce 
qui termine la preuve du lemme.
\end{preuve}

\begin{prop}
Le morphisme $\beta^\star$ envoie $\Fil^1 \Mod(\G)^\vee$ sur $\Fil^1 
\Mod(\G^\vee)$.
\end{prop}

\begin{preuve}
Soit $f \in \Fil^1 \Mod(\G)^\vee$. Par le lemme précédent, il existe un
morphisme de faisceaux $\tilde g$ faisant commuter le diagramme suivant :
$$\xymatrix @C=50pt {
\hhom (\G_n, \O_n^\cris) \ar[r]^-{\tilde f = \gamma_2^{-1}(f)} & 
\O_n^\cris \\
\hhom (\G_n, \J_n^\cris) \ar[r]^-{\tilde g} \ar@{^(->}[u] & \J_n^\cris 
\ar@{^(->}[u] }$$
Par ailleurs, le morphisme $\log : 1 + \J_n^\cris \to \O_n^\cris$ prend 
ses valeurs dans $\J_n^\cris$ et donc $\beta_\syn$ se factorise de la 
façon suivante :
$$\xymatrix @C=50pt {
\G_n^\vee \ar[r]^-{\beta_\syn} \ar@{.>}[rd] & \hhom (\G_n, \O_n^\cris) 
\\ & \hhom (\G_n, \J_n^\cris) \ar@{^(->}[u] }$$
En concaténant les deux diagrammes précédents, on remarque immédiatement 
que la composée $\tilde f \circ \beta_\syn$ tombe dans $\J_n^\cris$. 
Ainsi $\beta_\syn^\star (\tilde f) \in \hhom(\G_n, \J_n^\cris)$ ce
qui implique $\beta^\star (f) \in \Fil^1 \Mod(\G^\vee)$ comme annoncé.
\end{preuve}

\subsubsection{Compatibilité à $\phi_1$}

Dans ce paragraphe, on prouve que $\beta^\star$ est compatible à 
$\phi_1$. Pour cela, on introduit $\S_n$ le sous-faisceau de 
$\J_n^\cris$ défini comme le noyau du morphisme $\phi_1 - \id : 
\J_n^\cris \to \O_n^\cris$. On commence par montrer un lemme concernant 
ce faisceau :

\begin{lemme}
\label{5:lem:betasn}
Le morphisme $\beta_\syn$ se factorise par $\G_n^\vee \to \hhom(\G_n, 
\S_n) \toinj \hhom(\G_n, \O_n^\cris).$
\end{lemme}

\begin{preuve}
Il s'agit d'un calcul local pour la topologie syntomique. On reprend les 
notations ($A^\infty$, \emph{etc}.) du paragraphe 
\ref{5:subsec:synlocal}. En particulier, on dispose d'une suite 
exacte :
$$\xymatrix {
0 \ar[r] & 1 + \J_n^\cris(A^\infty) \ar[r] & W_n^{\cris, \DP} 
(A^\infty)^\star \ar[r]^-{s} & (A^\infty)^\star \ar[r] & 0 }.$$
D'après la description faite au début du paragraphe 
\ref{5:subsec:betasyn}, il suffit de montrer que si $x \in 
(A^\infty)^\star$ vérifie $x^{p^n} = 1$ et si $y \in W_n^{\cris, \DP} 
(A^\infty)^\star$ désigne un antécédent (que l'on a le droit de choisir)  
de $x$ alors $\phi_1(\log(y^{p^n})) = \log(y^{p^n})$. Considérons un tel 
$x$.

Il est immédiat de vérifier que $x$ admet un antécédent de la forme $y = 
(y_0, \ldots, y_{n-1}) \in W_n(A^\infty/p A^\infty)$. On note $z = 
y^{p^n} = (y_0^{p^n}, 0, \ldots, 0)$ et $\hat z = (y_0^{p^{n-1}}, 0, 
\ldots, 0) \in W_{n+1}(A^\infty/p A^\infty)$. Par définition, $\phi_1 
(\log(y^{p^n}))$ est la réduction modulo $p^n$ de $\frac 1 p 
\phi(\log \hat z)$. Mais :
$$\phi(\log \hat z ) = \log(\phi(\hat z)) = \log(\hat z^p) = p 
\, \log \hat z$$
l'égalité $\phi(\hat z) = \hat z^p$ étant vérifiée car $\hat z$ est un 
représentant de Teichmüller. On en déduit bien $\phi_1(\log(y^{p^n})) = 
\log(y^{p^n})$ comme voulu.
\end{preuve}

\begin{lemme}
\label{5:lem:tildeh}
Soient $f \in \Fil^1 \Mod(\G)^\vee$ et $h = \phi_1 (f)$. Alors, le 
diagramme de faisceaux sur $(T_{n+1})_\syn$ suivant :
$$ \xymatrix @C=50pt {
\hhom (\G_n, \J_n^\cris) \ar[r]^-{\gamma_2^{-1}(f)} \ar[d]_-{\phi_1} & 
\J_n^\cris \ar[d]^-{\phi_1} \\
\hhom (\G_n, \O_n^\cris) \ar[r]^-{\gamma_2^{-1}(h)} & \O_n^\cris } $$
est commutatif.
\end{lemme}

\begin{preuve}
Notons $\tilde f = \gamma_2^{-1}(f)$ et $\tilde h = \gamma_2^{-1}(h)$.
On fait à nouveau un calcul local : on reprend les notations 
($A^\infty$, \emph{etc}.) du paragraphe \ref{5:subsec:synlocal}.
Comme $\hhom (\G_n, \O_{T_n/E_n})$ est un cristal, on a :
$$w_\star \hhom (\G_n, \O_{T_n/E_n}) \simeq \Mod(\G) \otimes_{S_n} 
\O_n^\cris$$
et le morphisme $\tilde f$ (resp. $\tilde h$) s'écrit sur $A^\infty$ :
$$ \tilde f_{A^\infty} : \xymatrix @C=30pt {
\hhom (\G_n, \O_n^\cris) (A^\infty) \ar[r]^-{\eta} &
\Mod(\G) \otimes_{S_n} \O_n^\cris (A^\infty) \ar[r]^-{f \otimes 
\id} & \O_n^\cris(A^\infty) } $$
$$\text{(resp. } \tilde h_{A^\infty} : \xymatrix @C=30pt {
\hhom (\G_n, \O_n^\cris) (A^\infty) \ar[r]^-{\eta} &
\Mod(\G) \otimes_{S_n} \O_n^\cris (A^\infty) \ar[r]^-{h \otimes 
\id} & \O_n^\cris(A^\infty) } \text{)}.$$
Notons :
$$K = \Fil^1 \Mod(\G) \otimes_{S_n} \O_n^\cris (A^\infty) + \Mod(\G) 
\otimes_{S_n} \J_n^\cris (A^\infty) \subset \Mod(\G) \otimes_{S_n} 
\O_n^\cris (A^\infty)$$
et montrons que $\eta$ envoie $\hhom (\G_n, \J_n^\cris) (A^\infty)$ dans 
$K$. En reprenant  les notations de la preuve du lemme 
\ref{5:lem:gammafil}, on a un diagramme commutatif :
$$\xymatrix {
0 \ar[r] & \hhom(\G_n, \J_n^\cris)(A^\infty) \ar[r] & 
\hhom(\G_n, \O_n^\cris)(A^\infty) \ar[r] \ar[d]_-{\eta} & 
\hhom(\G_n, \O_n)(A^\infty) \ar[d] \\
& & \Mod(\G) \otimes_{S_n} \O_n^\cris(A^\infty) \ar[r]^-{\pr_\star 
\otimes s} & \calA_n^c \otimes_{\O_K/p^n} A^\infty }$$
où $s$ est la flèche définie en \ref{5:subsec:synlocal} et où $\pr_\star 
: \Mod(\G) \to \calA_n^c$ était noté $\pr_{\star(T_n,E_n)}$ dans le 
lemme \ref{5:lem:gammafil}. Il suffit donc de montrer que le noyau de 
$\pr_\star \otimes s$ est inclus dans $K$. Pour cela, on rappelle que
l'on a une suite exacte :
$$\xymatrix {
0 \ar[r] & \Fil^1 \Mod(\G) \ar[r] & \Mod(\G) \ar[r] & \calA_n^c \ar[r] 
& 0}$$
la surjectivité résultant du corollaire 3.2.10 de \cite{breuil} (avec 
$\H(n)$ à la place de $\H(1)$). Elle fournit le diagramme suivant :
$$\small \xymatrix @C=20pt @R=20pt {
& & \Mod(\G) \otimes_{S_n} \J_n^\cris(A^\infty) \ar[r] \ar[d] &
\calA_n^c \otimes_{S_n} \J_n^\cris(A^\infty) \ar[r] \ar[d] & 0 \\
0 \ar[r] & \Fil^1 \Mod(\G) \otimes_{S_n} \O_n^\cris(A^\infty) \ar[r]
& \Mod(\G) \otimes_{S_n} \O_n^\cris(A^\infty) \ar[r] \ar[d] 
\ar@{.>}[rd]^{\eta} & \calA_n^c \otimes_{S_n} \O_n^\cris(A^\infty) 
\ar[r] \ar[d] & 0 \\
& & \Mod(\G) \otimes_{S_n} A^\infty \ar[r] \ar[d] &
\calA_n^c \otimes_{S_n} A^\infty \ar[r] \ar[d] & 0 \\
& & 0 & 0 } $$
où toutes les lignes et colonnes sont exactes, l'exactitude de la ligne 
centrale résultant de la platitude de $\O_n^\cris(A^\infty)$ sur $S_n$. 
Une chasse au diagramme donne alors le résultat : le noyau de $\eta$ est 
inclus dans $K$.

On voit que $f \otimes \id$ envoie $K$ sur $\J_n^\cris (A^\infty)$ et 
que la restriction $\tilde f : \hhom(\G_n, \J_n^\cris) \to \J_n^\cris$ 
s'écrit sur $A^\infty$ de la façon suivante :
$$ \tilde f_{A^\infty} : \xymatrix @C=30pt {
\hhom (\G_n, \J_n^\cris) (A^\infty) \ar[r]^-{\eta} &
K \ar[r]^-{f \otimes \id} & \J_n^\cris(A^\infty) }. $$
On considère pour terminer le diagramme suivant :
$$ \xymatrix @C=50pt {
\hhom (\G_n, \J_n^\cris) (A^\infty) \ar[r]^-{\eta} \ar[d]_-{\phi_1} &
K \ar[r]^-{f \otimes \id} \ar[d]_-{\phi_1 \otimes \phi + \phi \otimes 
\phi_1} & \J_n^\cris(A^\infty) \ar[d]^-{\phi_1} \\
\hhom (\G_n, \O_n^\cris) (A^\infty) \ar[r]^-{\eta} &
\Mod(\G) \otimes_{S_n} \O_n^\cris (A^\infty) \ar[r]^-{h \otimes 
\id} & \O_n^\cris(A^\infty) }$$
où $\phi : \Mod(\hat \G) \to \Mod(\hat \G)$ est défini par la formule 
$\phi(x) = \frac 1 c \phi_1 (E(u) x))$. On vérifie que les deux carrés
commutent. On en déduit que tout le diagramme commute, ce qui démontre 
le lemme.
\end{preuve}

\bigskip

\noindent
{\it Remarque.} En examinant la preuve précédente, on constate qu'elle 
fournit une autre démonstration (pas très éloignée toutefois) du
lemme \ref{5:lem:gammafil}.

\medskip

\begin{prop}
Le morphisme $\beta^\star$ est compatible à $\phi_1$.
\end{prop}

\begin{preuve}
Pour cette preuve, on note encore $\beta_\syn$ le morphisme
$\G_n^\vee \to \hhom(\G_n, \J_n^\cris)$ induit par $\beta_\syn$.

Soit $f \in \Fil^1 \Mod(\G)^\vee$. Notons $\tilde f = \gamma_2^{-1}(f)$. 
D'après le lemme \ref{5:lem:gammafil}, il existe un morphisme $\tilde 
g : \hhom(\G_n, \J_n^\cris) \to \J_n^\cris$ faisant commuter le 
diagramme suivant :
$$\xymatrix @C=50pt {
\hhom (\G_n, \O_n^\cris) \ar[r]^-{\tilde f} & \O_n^\cris \\
\hhom (\G_n, \J_n^\cris) \ar[r]^-{\tilde g} \ar@{^(->}[u] & \J_n^\cris 
\ar@{^(->}[u] }$$
Le morphisme $\phi_1 (\beta^\star f)$ est alors défini comme la composée 
$\phi_1 \circ \tilde g \circ \beta_\syn$. On considère le diagramme 
suivant :
$$\xymatrix @C=60pt {
\G_n^\vee \ar[r]^-{\beta_\syn} \ar@{=}[d] & \hhom (\G_n, \J_n^\cris) 
\ar[r]^-{\tilde g} \ar[d]^-{\phi_1} & \J_n^\cris \ar[d]^-{\phi_1} \\
\G_n^\vee \ar[r]^-{\beta_\syn} & \hhom (\G_n, \O_n^\cris) 
\ar[r]^-{\gamma_2^{-1} (\phi_1 (f))} & \O_n^\cris }$$
Le carré de gauche commute (lemme \ref{5:lem:betasn}), ainsi que celui 
de droite (lemme \ref{5:lem:tildeh}). On en déduit que tout le diagramme 
commute et donc que :
$$\beta^\star (\phi_1 (f)) = \beta_\syn^\star \circ \gamma_2^{-1}
(\phi_1 (f)) = \gamma_2^{-1} (\phi_1 (f)) \circ \beta_\syn =
\phi_1 \circ \tilde g \circ \beta_\syn = \phi_1 (\beta^\star f)$$
ce qui implique la proposition.
\end{preuve}

\subsection{Cas général}
Le but de ce dernier paragraphe est de prouver le théorème suivant :

\begin{theo}
Soit $\G$ un $\O_K$-groupe. Alors on a un isomorphisme canonique et 
fonctoriel :
$$\Mod(\G)^\vee \simeq \Mod(\G^\vee)$$
dans la catégorie $\Mods$.
\end{theo}

\begin{preuve}
Bien sûr, l'isomorphisme dont il est question est $\beta^\star$ défini 
en \ref{5:sec:morcomp}.

\medskip

On montre dans un premier temps que le morphisme $\beta^\star$ est 
compatible à $\Fil^1$ et $\phi_1$. D'après ce que l'on a fait 
précédemment c'est vrai si $\G$ est de la forme $\H(n)$ pour un groupe
$p$-divisible $\H$. Dans le cas général, il existe un épimorphisme 
$\H(n) \to \G$ pour un 
certain groupe $p$-divisible $\H$ et un certain entier $n$. D'après la 
proposition \ref{5:prop:dualmr} et la proposition 4.2.1.5 de 
\cite{breuil}, on a un diagramme commutatif :
$$\xymatrix @C=50pt {
\Mod(\H(n))^\vee \ar[r]^-{u} \ar[d]_-{\beta^\star_{\H(n)}} &  
\Mod(\G)^\vee \ar[d]^-{\beta^\star_{\G}} \\
\Mod(\H(n)^\vee) \ar[r]^-{v} & \Mod(\G^\vee) }$$
et la flèche induite $u : \Fil^1 \Mod(\H(n))^\vee \to \Fil^1 
\Mod(\G)^\vee$ est surjective. Soient $x \in \Fil^1 \Mod(\G)^\vee$ et $y
\in \Fil^1 \Mod(\H(n))^\vee$ un antécédent par $u$ de $x$. On a 
$\beta^\star_{\G} (x) = v \circ \beta^\star_{\H(n)} (y)$ d'où, puisque 
$v$ et $\beta^\star_{\H(n)}$ respectent le $\Fil^1$, il vient 
$\beta^\star_{\G} (x) \in \Fil^1 \Mod(\G^\vee)$ et la compatibilité 
recherchée.

\medskip

Passons à la compatibilité avec $\phi_1$. On considère cette fois-ci 
un monomorphisme $\G \to \H(n)$, qui donne naissance au cube commutatif 
suivant :
\[\xymatrix @C=30pt @R=20pt {
& \Mod(\G)^\vee \ar@{^(->}[rr] \ar[dd] & & \Mod(\H(n))^\vee
\ar[dd]^{\beta^\star_{\H(n)}} \\
\Fil^1 \Mod(\G)^\vee \ar@{^(->}[rr] \ar[dd]_{\beta^\star_{\G}} 
\ar[ur]^{\phi_1} & \ar@{ }[d]_{\beta^\star_{\G}} & \Fil^1 
\Mod(\H(n))^\vee  \ar[dd] \ar[ur]^{\phi_1} \\
& \Mod(\G^\vee) \ar@{^(->}[rr] & \ar@{ }[d]_{\beta^\star_{\H(n)}} & 
\Mod(\H(n)^\vee) \\
\Fil^1 \Mod(\G^\vee) \ar@{^(->}[rr] \ar[ur]^{\phi_1} & & \Fil^1 
\Mod(\H(n)^\vee) \ar[ur]^{\phi_1} }\]
Les flèches horizontales sont toutes injectives et toutes les faces, 
sauf \emph{a priori} celle de gauche, commutent. Une chasse au diagramme 
permet de prouver que la face de gauche est aussi commutative (on 
utilise l'injectivité de $\Mod(\G)^\vee \to \Mod(\H(n))^\vee$) et donc 
de conclure.

\medskip

Finalement, il ne reste plus qu'à prouver que $\beta^\star : 
\Fil^1 \Mod(\G)^\vee \to \Fil^1 \Mod(\G^\vee)$ est surjectif. Or, par
le lemme (facile) 4.2.14 de \cite{breuil}, cela est automatique lorsque
$\G$ est tué par $p$. On conclut par un dévissage aisé laissé au 
lecteur.
\end{preuve}

\bigskip

Par passage à la limite, on en déduit le théorème suivant :

\begin{theo}
Soit $\G$ un groupe $p$-divisible sur $\O_K$. Alors, si $G^\vee$ désigne 
le dual de Cartier de $\G$, on a un isomorphisme canonique et  
fonctoriel :
$$\Mod(\G)^\vee \simeq \Mod(\G^\vee)$$
dans la catégorie des modules fortement divisibles.
\end{theo}

\bibliography{cartier}

\renewcommand{\refname}{Bibliographie}
 
\providecommand{\bysame}{\leavevmode\hbox to3em{\hrulefill}\thinspace}

\begin{thebibliography}{\hspace{1.3cm}}
\bibitem[BBM82]{bbm}
 P. Berthelot, L. Breen et W. Messing,
 {\it Théorie de Dieudonné cristalline II},
 Lecture notes in math. 
 {\bf 930},
 Springer-Verlag 
 (1982)

\bibitem[Bre98]{breuil-duke}
 C. Breuil,
 {\it Cohomologie étale de $p$-torsion et cohomologie cristalline en réduction semi-stable},
 Duke mathematical journal 
 {\bf 95}
 (1998), 
 523--620

\bibitem[Bre99]{breuil-invent}
 \bysame ,
 {\it Représentation semi-stables et modules fortement divisibles},
 Invent. math. 
 {\bf 136}
 (1999), 
 89--122

\bibitem[Bre00]{breuil}
 \bysame ,
 {\it Groupes $p$-divisibles, groupes finis et modules filtrés},
 Annals of Mathematics 
 {\bf 152}
 (2000), 
 489--549

\bibitem[Car05]{caruso-these}
 X. Caruso,
 {\it Conjecture de l'inertie modérée de Serre},
 thèse 
 (2005)

\end{thebibliography}
\bibliographystyle{amsalpha}

\end{document}